\documentstyle[12pt]{article}
\newcommand{\const}{\mathop{\rm const}\limits}

\newcommand{\supp}{\mathop{\rm supp}\limits}

\newcommand{\vraisup}{\mathop{\rm vraisup}\limits}

\begin{document}

\begin{center}

{\bf  SHARP VALUE FOR THE NORM OF COMPOSITION} \\

\vspace{4mm}

{\bf AND MULTIPLICATIVE OPERATORS BETWEEN } \\

\vspace{4mm}

{\bf TWO DIFFERENT GRAND LEBESGUE SPACES} \par

\vspace{4mm}

{\bf E. Ostrovsky}\\

e-mail: eugostrovsky@list.ru \\

\vspace{4mm}

{\bf L. Sirota}\\

e-mail: sirota3@bezeqint.net \\

\vspace{3mm}

Department of Mathematics and Statistics, Bar-Ilan University, 59200, \\
 Ramat Gan, Israel. \\

\vspace{3mm}

{\bf Abstract.} \\

\vspace{3mm}

 \ We calculate in this paper the norm of composition, multiplicative and product operator, generated by multiplicative and
 measurable argument transformation between two different ordinary Lebesgue-Riesz and Grand Lebesgue spaces. \par

 \ We set ourselves the aim to obtain the exact expression for the norm of the considered operators by means of
building of appropriate examples. \par

\vspace{2mm}

\end{center}

\vspace{2mm}


\vspace{3mm}


\section{Introduction}

\vspace{2mm}

 \  Let $  (X = \{x\}, M, \mu)  $ and $ (Y = \{y\}, N, \nu )  $  be two measurable spaces equipped with a non-zero sigma-finite
measures $  \mu $ and $  \nu $ correspondingly. Denote by $ N_0  = N_0(Y) $   ( and correspondingly $ \ M_0 = M_0(X) \ ) \ $
the linear set of all numerical measurable functions $ \ f: Y \to R. \ $
The case $ X = Y $ can not be excluded. \par

 \  Let also $  \xi = \xi(x)  $ be measurable function from the set $  X   $ to $  Y: \ \xi: X \to Y.  $ \par

 \vspace{3mm}

 {\bf Definition 1.1.} \par

\vspace{3mm}

 \  The linear operator  $  U_{\xi} [f] = U_{\xi} [f](x) $ defined may be on some Banach invariant
subspace $  B_1 $ of the space  $  M_0(X) \ $
to another one, in general case, subspace $  B_2 $  of the space $ N_0(Y) $  by the formula

 $$
 U_{\xi} [f](x) = f(\xi(x))  \eqno(1.1)
 $$
is said to be {\it  composition operator }  generated by $  \xi(\cdot). $ \par

\vspace{3mm}

 \ Many  important properties on these operators acting on different spaces  $ B_1 $  with values in $  B_2: $ Lebesgue, Lorentz,
mainly  Orlicz  spaces etc., namely: boundedness, compactness, the exact values of norm  are investigated, e.g. in
 \cite{Arora1}, \cite{Cui1}, \cite{Douglas1}, \cite{Estaremi1}, \cite{Estaremi2}, \cite{Estaremi101}, \cite{Kumar1},
\cite{Montes-Rodrguez1}, \cite{Singh1}, \cite{Gupta1}, \cite{Shapiro1}, \cite{Takagi1} - \cite{Takagi4}, \cite{Zheng1}. \par
 \ Applications other than those mentioned above appears in ergodic theory, see, e.g., in
 \cite{Anosov1},  \cite{Bedford1}, \cite{Petersen1}, \cite{Walters1} etc. \par

 \ We study first of all the action of these operators  between the classical Lebesgue-Riesz spaces $ L_p(X,\mu) = L_{p,X} $ consisting
on all the measurable functions with finite norm

$$
|\xi|L_p(X,\mu) = |\xi| L_{p,X} = | \ \xi \ |_p := \left[\int_X  |\xi(x)|^p \ \mu(dx) \right]^{1/p}, \ p = \const \ge 1, \eqno(1.2a)
$$
 and correspondingly

$$
|g|L_q(Y,\nu) = |g| L_{q,Y} = | g|_q := \left[\int_Y  |g(y)|^q \ \nu(dy) \right]^{1/q}, \ q = \const \ge 1. \eqno(1.2b)
$$

$$
|f|L_r(Y,\nu) = |f| L_{r,Y} = | f|_r := \left[\int_Y  |f(y)|^r \ \nu(dy) \right]^{1/r}, \ r = \const \ge 1, \eqno(1.2c)
$$

 \ As ordinary,

$$
|\xi|L_{\infty}(X,\mu) = |\xi| L_{\infty,X} = | \ \xi \ |_{\infty} := \vraisup_{x \in X} |\xi(x)|, \eqno(1.2d)
$$
and analogously may be defined the norm  $  |g|L_{\infty}(Y,\nu) = |g| L_{\infty,Y} = | g |_{\infty}.  $ \par

 \vspace{4mm}

 \ {\bf  Our purpose in this short article is calculation  of the exact value of the norm for these operators in ordinary Lebesgue-Riesz
 as well as in Grand Lebesgue spaces (GLS):}

\vspace{3mm}

$$
|U_{\xi}|_{q \to p} =
R(q,p) = R(q,p; Y,X) \stackrel{def}{=} \sup_{0 \ne f \in L_q(Y,\nu)} \left[ \frac{|U_{\xi}[f]|_{p,X}}{|f|_{q,Y}} \right]. \eqno(1.3)
$$

\vspace{3mm}

\ {\bf We consider further at the same problems for some another operators: multiplicative and product of these operators.}\par

\vspace{3mm}

 \ The case of the so-called Grand Lebesgue Spaces (GLS) will be also considered further. \par

 \ Some estimations of the norm of these operators acting between Orlicz spaces are obtained in the articles
\cite{Estaremi1}, \cite{Estaremi2}, \cite{Estaremi101}. A particular case $ \xi(x) = 1/x, \ X = (0,\infty) $ is considered
in a famous monograph \cite{Okikiolu1},  pp. 220 - 221. See also \cite{Dunford1}, chapter n 7, pp. 660-666, where is
considered the case  $ \ X = Y \ $ and $ \ q = p.  $ \par

 \ The classical (kernel) integral operators, including singular, acting
  in Lebesgue-Riesz, Orlicz, Grand Lebesgue spaces are investigated in \cite{Okikiolu1},
\cite{Rao1}, p. 198-220, \cite{Ostrovsky8} etc.

\vspace{3mm}

 \ Another operators acting in these spaces: Hardy, Riesz, Fourier, maximal, potential etc. are investigated, e.g.  in
 \cite{Liflyand1}, \cite{Ostrovsky100}, \cite{Ostrovsky101}, \cite{Ostrovsky3}, \cite{Ostrovsky4}. \par

 \ Other notations. Denote by  $ F(\cdot) = F_{\xi}(\cdot) $ the {\it distribution} of the (measurable) function $ \xi: $

$$
F(A) = F_{\xi}(A) = \mu(x: \ \xi(x) \in A ), \ A \in N; \eqno(1.4)
$$
then  $ F_{\xi}(\cdot)  $ is sigma-additive and sigma-finite measure in $  N. $ \par
 \ Introduce also the Radon-Nikodym  derivative $ z = z(y) = dF /d\nu $ of the measure $  F_{\xi} $  relative the source
measure $ \nu $ on the $  Y, $ i.e. such that for arbitrary measurable function $ h: Y \to R  $

$$
\int_X h(\xi(x)) \ \mu(dx) = \int_Y h(y) \ z(y) \ \nu(dy). \eqno(1.5)
$$
 \ If the measure $ F_{\xi} $ is not absolutely  continue relative the measure $ \nu, $ we define formally $ z(y) = + \infty. $\par

 \ In particular,

$$
 | \ |U_{\xi}[f] \ |_p^p  =   \int_X |f(\xi(x))|^p \ \mu(dx) =  \int_Y |f(y)|^p \ z(y) \ \nu(dy). \eqno(1.6)
$$

\bigskip

\section{ Main result. The case of ordinary Lebesgue-Riesz spaces. }

\vspace{3mm}

 \ Define also the following important functional

$$
K_{z}(p,q) \stackrel{def}{=} | \ z \ |_{q/(q-p), Y}^{1/p} =
\left[\int_Y z(y)^{ q/(q-p)} \ \nu(dy)  \right]^{1/p - 1/q}, \eqno(2.1)
$$
in the case when $  \ q = \const > p = \const \ge 1 $ and

$$
K_{z}(p,p) \stackrel{def}{=} | \ z \ |_{\infty,Y}^{1/p} =
[ \ \vraisup_{y \in Y} z(y) \ ]^{1/p},  \eqno(2.1a)
$$
if $  p = \const \ge 1. $\par

\vspace{3mm}

\ {\bf Remark 2.1.} We does not exclude the case when the integral in (2.1) divergent; then evidently $ K_{z}(p,q) = + \infty. $ \par

\vspace{3mm}

 \ {\bf Theorem 2.1.}

$$
|U_{\xi}|_{q \to p} =: R(q,p; Y,X) = K_{z}(p,q). \eqno(2.2)
$$

\vspace{3mm}

 \ {\bf Proof. Upper bound.}  It is sufficient to consider the case $  1 \le p < q < \infty. $
Suppose $ K_{r,z}(p,q) < \infty; \ $ the opposite case will be considered further. \par

 \ We apply the H\"older's inequality to the right-hand side of inequality (1.6):

$$
 | \ |U_{\xi}[f] \ |_p^p  \le  \left[   \int_Y |f(y)|^{ \alpha p} \ \nu(dy) \right]^{1/\alpha} \cdot
  \left[ \int_Y \ z^{\gamma}(y) \ \nu(dy) \right]^{1/\gamma} =
$$

$$
| \ f \ |^p_{\alpha p, Y} \cdot | \ z \ |_{\gamma, Y},\eqno(2.3)
$$
where $ \alpha = \const > 1 $ and $ \gamma $ is its conjugate number $ \gamma = \alpha/(\alpha - 1). $ \par

 \ If we choose in (2.3) $  \alpha = q/p > 1, $ then

$$
\gamma = \frac{q}{q  - p}.
$$

 \ Substituting into (2.3), we get  after simple calculations to the estimate

$$
|U_{\xi}[f]|_{p,X} \le  K_{z}(p,q) \cdot |f|_{q,Y}. \eqno(2.4)
$$

\vspace{3mm}

 \ It remains to  prove that the condition $  F_{\xi}(\cdot) << \nu(\cdot) $ is necessary for the inequality

$$
|U_{\xi}|_{q \to p} = R(q,p) = R(q,p; Y,X) < \infty
$$
 at last for {\it some } pairs of positive numbers $ p $ and $  q. $   Suppose there exists a finite constant $ W = W(p,q) $
such that

$$
\left[ \ \int_Y |f(y)|^p \ F_{\xi}(dy) \ \right]^{1/p} \le W(p,q) \  \left[ \ \int_Y |f(y)|^q \ \nu(dy) \ \right]^{1/q}  \eqno(2.5)
$$
for arbitrary function $ f (\cdot) \in L_q(Y,\nu).  $  \par

 \ We can and will assume as the capacity of the value $  W(p,q) $ for all the positive values $  \ p,q \  $
 its minimal value, indeed

$$
W(p,q) := \sup_{0 < |f|_{q,Y} < \infty } \left\{ \frac{ [ \ \int_Y |f(y)|^p \ F_{\xi}(dy) \ ]^{1/p}}{|f|_{q,Y}}  \right\}.
$$

 \ We deduce choosing $  f(y) = I_B(y), $ where  $ I(\cdot) $ is ordinary indicator function of the measurable set
$ B \in N $ and $  0 \le \nu(B) < \infty: $

$$
\left[ \ F_{\xi}(B) \ \right]^{1/p} \le W(p,q) \ \left[ \ \nu(B) \ \right]^{1/q}.
$$

 \ Therefore,  every time when $ \nu(B) = 0,  $ then  right here $  F_{\xi}(B) = 0. $ Thus, the distribution measure $ F_{\xi}(\cdot) $
is absolutely continuous relative the measure  $ \nu(\cdot): \hspace{5mm} F_{\xi}(\cdot) << \nu(\cdot). $\par

 \ Note in addition that in the case $ q = p $ the correspondent Radon-Nikodym derivative is bounded:

$$
\frac{dF_{\xi}}{d \nu}(y) \le W(p,p)
$$
and following in this case

$$
|| \ U_{\xi} \ ||(L_{p,Y} \to L_{p,X}) = W(p,p).
$$

 \ This fact was proved first in the  famous  article of Cui, Hudzik, Kumar, Maligranda \ \cite{Cui1}.\par

 \ We will see further that in the general case $  1 \le p < q < \infty $ the Radon-Nikodym  derivative $ dF_{\xi}/d \nu $
may be unbounded. \par

\vspace{4mm}

{\bf Proof. Lower bound.} Let us choose as a capacity of the trial function

$$
f(y) := z^{\beta}(y), \ \beta = \const.  \eqno(2.6)
$$
 \ We have

$$
|f|_{q,Y} =  \left[ \int_Y  z^{\beta q}(y) \ \nu(dy)  \right]^{1/q} = |  \ z \ |^{\beta}_{\beta q, Y},
$$

$$
| \ u_{\xi}[f] \ |_{p,X} = \int_X |f(y)|^p \ z(y) \ \nu(dy) = \int_Y z^{\beta p + 1} \ \nu(dy) =
|  \ z \ |_{\beta p + 1, Y}^{\beta p + 1};
$$

$$
Q[f] := \frac{| \ u_{\xi}[f] \ |_{p,X}}{|f|_{q,Y}} = \frac{|  \ z \ |_{\beta p + 1, Y}^{\beta p + 1}}{|  \ z \ |^{\beta}_{\beta q, Y}}.
 \eqno(2.7)
$$

 \ If we choose the value $  \beta $ such that $ \beta p + 1 = \beta q, $ i.e. $ \beta := 1/( q - p ), $ then

$$
Q[f] = | \ z \ |^{1/p}_{ q/(q - p)} = K_{z}(p,q).  \eqno(2.8)
$$

 \ The case when $ K_z(p,q) = \infty $ may be investigated analogously. Namely,
let us choose the {\it sequence} of truncated  trial functions

$$
f_n(y) := z^{\beta}(y) \ I( z (y) \le n), \ n = 1,2,\ldots;  \ \beta = \const = 1/(q - p),
$$
where $ I(\cdot) $ is again indicator function. We have as before as $  n \to \infty $

$$
Q[f_n] = \frac{| \ u_{\xi}[f_n] \ |_{p,X}}{|f_n|_{q,Y}}  =
 \frac{|  \ z \ I(z(y) \le n) \ |_{\beta p + 1, Y}^{\beta p + 1}}{|  \ z \ I(z(y) \le n) \ |^{\beta}_{\beta q, Y}}  =
$$

$$
= | \ z \ I(z(y) \le n \ |^{1/p}_{ q/(q - p)} \to | \ z \ |^{1/p}_{ q/(q - p)} =  K_{z}(p,q) = \infty.  \eqno(2.9)
$$

 \ Thus, and in this case

$$
|U_{\xi}|_{q \to p} = \infty = K_z(p,q).
$$

 \vspace{3mm}

 \ This completes the proof of theorem 2.1. \par

\vspace{3mm}

\ {\bf  Example 2.1.} Let $ X = Y = [0,1] $ with ordinary Lebesgue measure and let $  r = \const > 0. $  We exclude also the trivial
case  $  r = 1. $\par
 \ Define the linear composition operator of the form

$$
U_{(r)}[f](x) = f(x^r). \eqno(2.10)
$$

 \ The correspondent Radon-Nikodym derivative $  \ z(y) = z^{(r)}(y) \ $ is equal here

$$
z(y) = z^{(r)}(y) = r^{-1} \ y^{ 1/r - 1 }, \ y > 0.
$$

 \ The expression for the "constant" $ K_z(p,q) = K_z^{(r)}(p,q) $ has here a form

$$
K_z^{(r)}(p,q) = r^{-1/p} \left( \ \frac{q - p}{q - pr} \ \right)^{1/p \ - 1/q} \stackrel{def}{=} S_r(p,q).  \eqno(2.11)
$$

 \ We  apply the proposition of theorem 2.1

$$
|| \ U_{(r)} \ ||(L_q[0,1] \to L_p[0,1]) = S_r(p,q). \eqno(2.12)
$$

 \ Let us  investigate in detail the  variable $ \ S_r(p,q). \ $
 We can distinguish the following two cases: a case A) $ \ 0 < r < 1, \ $ (a "good case") and a case B) $ \ r > 1 \ $ ("a hard  case").\par

 \ In the first case A) the "transfer" function $ K_{z}^{(r)}(p,q) $ is bounded relative the variable $  q  $ in the {\it closed} domain
 $ q \in [p, \infty]: $

$$
\sup_{ q \in [p, \infty] } K_{z}^{(r)}(p,q) \le r^{-1/p} \eqno(2.13)
$$
and wherein

$$
\lim_{q \to p+0}  K_{z}^{(r)}(p,q) = r^{-1/p}, \eqno(2.14)
$$

$$
K_{z}^{(r)}(p,p) = r^{-1/p}  =\left[ \ \vraisup_{y \in (0,1)} z^{(r)}(y) \ \right]^{1/p}. \eqno(2.14a)
$$

\vspace{3mm}

 \ Let us now consider the opposite case $  \ r > 1. \ $ Then the "transfer" function $ K_{z}^{(r)}(p,q) $ is finite only
in the {\it open} interval $ \ q  > pr \ $  and we have as $  q \to pr + 0 $

$$
K_{z}^{(r)}(p,q) \sim r^{-1/p} \ \left[ \frac{p(r-1)}{q-pr} \right]^{(r-1)/pr }. \eqno(2.15)
$$
 \ Thus, in this case the "transfer"  function $ \ K_{z}^{(r)}(p,q) \ $ is really unbounded.\par

\bigskip

\section{Multiplicative operators.}

\vspace{3mm}

 \ Define the so-called {\it multiplicative } operator $  V  = V_g $ acting from the space $  \ N_0(Y) $ into itself

$$
V_g[f] (y) \stackrel{def}{=} g(y) \cdot f(y),\eqno(3.1)
$$
where a factor $ g = g(y) $ is measurable function ("weight"):  $ g: Y \to R. $ \par

 \ The acting of these operators between Orlicz's spaces, satisfying as a rule the $ \ \Delta_2 \ $ condition, was investigated
in many works, see e.g. \cite{Abrahamse1}, \cite{Gupta1}, \cite{Musielak1}; see also a recent article \cite{Giri1} and reference
therein.\par

 \ Denote also

$$
Q_g(p,q) := | \ g \ |_{pq/(q - p)} = | \ g \ |_{pq/(q - p),Y}.
$$

\vspace{3mm}

\ {\bf Proposition 3.1.}

\vspace{3mm}

$$
|| \ V_g \ ||{q \to p} = Q_g(p,q) = |g|_{pq/(q - p)}, \ q > p, \eqno(3.2)
$$

\vspace{3mm}

and $ || \ V_g \ ||{q \to p} = \infty $ if $ q \le p $ or when $  g \notin L_{pq/(q-p)}. $ \par

\vspace{3mm}

\ {\bf  Proof } is quite analogous to ones in theorem 2.1; we omit some non-essential details. Note first of all
that again by virtue of H\"older's inequality

$$
| \ g \ f  \ |_p  \le | \ f \ |_q \cdot | \ g \ |_r,
$$
where

$$
\frac{1}{p} = \frac{1}{q} +  \frac{1}{r}, \hspace{4mm} p,q,r > 1.
$$
 \ Since $ r =pq/(q - p), \ q > p, $ we conclude

$$
|| \ V_g \ ||_{q \to p} \le |g|_{pq/(q - p)}, \ q > p.\eqno(3.3)
$$

 \vspace{3mm}

 \ It remains to obtain the lower estimate. Suppose without loss of generality $ g(x) \ge 0. $
One can choose the correspondent trial function  as follows

$$
f(y) := g^{\beta}(y), \ \beta = \const > 0.
$$
 Then

$$
| \ V_g[f] \ |_p^p = \int_Y g^{ (\beta + 1)p }(y) \ \nu(dy) = | \  g(\cdot) \ |_{(\beta + 1)p, Y}^{(\beta + 1)p};
$$

$$
| \ f \ |_{q,Y} = \left[ \int_Y g^{\beta q}(y) \ \nu(dy) \right]^{1/q} = | \ g \ |^{\beta}_{\beta q, Y};
$$

$$
\frac{| \ V_g[f] \ |_{p,Y} }{|f|_{q,Y}} = \frac{| \ g \ |_{(\beta + 1)p,Y}^{\beta + 1}}{| \ g \ |^{\beta}_{\beta q, Y}};
$$
and when we choose $ \beta = p/(q - p), $ we obtain what is desired:

$$
|| \ V_g \ ||_{q \to p} \ge |g|_{pq/(q - p)}, \ q > p.\eqno(3.4)
$$

\vspace{3mm}

{\bf Example 3.1.} Let again $  Y = (0,1) $ with ordinary Lebesgue measure and define the following multiplicative operator

$$
V^{(t)}[f](y) \stackrel{def}{=} y^{-t} \cdot f(y),   \ t = \const \in (0,1);
\eqno(3.5)
$$
the case $  t \le 0 $ is trivial for us. I.e. in this case

$$
g(y) = g^{(t)}(y) := y^{-t}.
$$

 \ We obtain after simple calculations

$$
|| \ V_{g^{(t)}} \ ||_{q \to p} =  Q_{g^{(t)}}(p,q) := \left( 1 - \frac{t p q}{ q - p} \right)^{1/q - 1/p} \eqno(3.6)
$$
in the case when

$$
 1 \le p < q, \hspace{7mm} \frac{t p q}{ q - p} < 1
$$
and

$$
 || \ V_{g^{(t)}} \ ||_{q \to p} = + \infty \eqno(3.7)
$$
otherwise. \par

\vspace{3mm}

 {\bf Remark 3.1.} Note that the case when the function $ \ g(\cdot) \  $  is bounded is trivial; in this case one can take $ \ p = q. $ \par

\vspace{3mm}

 {\bf Remark 3.2.} Many another examples may be found further in the seventh  section. \par

\bigskip

\section{Product operators.}

\vspace{3mm}

 \ Let us consider in this section the so-called "product" operator $ W_{g,\xi}[f] $ of the form

$$
W_{g,\xi}[f](x) = g(x) \cdot f(\xi(x)), \hspace{4mm} x \in X, \eqno(4.1)
$$
where as above $  \ \xi: X \to Y, \  f: Y \to R, \ g: X \to R  $ be measurable functions. \par
 \ On the other words, the operator $  W_{g,\xi}(\cdot)  $ is a product of two non commuting, in general case, operators
$  W_{g,\xi}(\cdot) = V_g(\cdot) \odot U_{\xi}(\cdot).$ \par

 \ These operators play a very important role in the study of linear isometries acting in the Lebesgue-Riesz spaces $  L_p, $ see
\cite{Okikiolu1}, p. 176-177. \par

 \ We retain all the assumptions and notations of all the foregoing sections. \par

\vspace{3mm}

{\bf A. General case.} \\

\vspace{2mm}

  \ Let $ f(\cdot) \in L_{q(1),Y}, \ q(1) > 1; $ and let $  q(2) < q(1). $ We conclude by virtue of
 theorem 2.1 $  U_{\xi}[f] \in L_{q(2),X}  $  and herewith

$$
| \ U_{\xi}[f] \ | L_{q(2),X}  \le K_z(q(1),q(2)) \ | \ f \ |L_{q(1),Y}.
$$

 \ We have further applying theorem 3.1  for the certain  value $  q(3) < q(2) $

$$
| \ W_{g,\xi}[f] \ | L_{q(3),X} \le Q_g(q(2), q(3)) \cdot K_z(q(1), q(2)) \cdot | \ f \ |L_{q(1),Y}. \eqno(4.2)
$$

 \ The estimate (4.2) may be issued as follows. Let us introduce the following function

$$
T(p,q) = T_{g,\xi}(p,q) \stackrel{def}{=}
\inf_{ l \in (p,q) } \left[ K_z(p,l) \ Q_g(l,q) \right], \ 1 \le p < q < \infty.
$$

\vspace{3mm}

\ {\bf Theorem 4.1.}  Let $  1 \le p < q < \infty. $  Then

$$
|| \ W_{g,\xi} \ ||(L_p(Y) \to L_q(X) ) \le  T_{g,\xi}(p,q)  =
 \inf_{ l \in (p,q) } \left[ K_z(p,l) \ Q_g(l,q) \right]. \eqno(4.3)
$$

\vspace{3mm}

 \ Note that at the same result may be obtained by means of H\"older's inequality applying to the right-hand side of
the relation (4.1). \par

\vspace{4mm}

{\bf B. Particular case.} \\

\vspace{2mm}

 \ Suppose here that the weight function (factor) $ \ g(x) \ $ has a form $  g(x) = h(\xi(x)), $ where
 as before $  \ \xi: X \to Y, \  h: Y \to R $ are measurable functions. On the other words, here

$$
W_{g,\xi}[f](x) = h(\xi(x)) \cdot f(\xi(x)). \eqno(4.4)
$$

 \ We deduce using umpteenth time H\"older's inequality for three multipliers

$$
| \ W_{g,\xi}[f] \ |_{p,X} \le | \ h \ |_{p \theta,Y} \ | \ f \ |_{q, Y} \ \left[ \ | \ z \ |_{\tau, Y} \ \right]^{1/p}, \eqno(4.5)
$$
where

$$
q > p,  \hspace{4mm} \theta, \tau > 1 \eqno(4.6a)
$$
and

$$
 \frac{1}{\theta} + \frac{1}{\tau} = \frac{q - p}{q}. \eqno(4.6b)
$$

\vspace{3mm}

 \ Eventually:  \par

\vspace{3mm}

\ {\bf Theorem 4.2.} We propose under formulated in this pilcrow definitions and conditions

$$
|| \ W_{g,\xi} \ || (L_{q,Y} \to L_{p,X}) =
\inf_{\theta, \tau} \left\{ \ | \ h \ |_{p \theta,Y} \cdot \left[ \ | \ z \ |_{\tau, Y} \ \right]^{1/p} \ \right\}, \eqno(4.7)
$$
where $ "\inf" $ in the inequality (4.7) is calculated over all  the variables satisfying  the relations (4.6a) and (4.6b). \par

\vspace{3mm}

{\bf C. Independent case.} \\

\vspace{2mm}

 \ Suppose in this subsection that both the functions  $  f(\xi(\cdot))  $ and $  g(\cdot) $ are independent in the theoretical
probability sense, i.e.

$$
\mu \{x: \ g(x) \in A_1, \ f(\xi(x) ) \in A_2  \}=   \mu \{x: \ g(x) \in A_1 \} \times
$$

$$
\mu \{x: \ f(\xi(x) ) \in A_2  \}, \ A_1, A_2 \in M.
$$
 \ Then

$$
| \ g \ f(\xi(\cdot)) \ |_{p,X}  = | \ g \ |_{p,X} \cdot | \ f(\xi(\cdot)) \ |_{p,X}. \eqno(4.8)
$$

 \ This statement may be issued in the considered here independent case as follows. \par

 \vspace{3mm}

 \ {\bf Theorem 4.3.} We propose under conditions of independent case by virtue of theorem 2.1

$$
|| \ W_{g,\xi} \ || (L_{q,Y} \to L_{p,X}) = |\ g \ |_p \cdot K_z(p,q). \eqno(4.9)
$$

\vspace{3mm}

 \ We note in conclusion that the lower bounds in the last two theorems are trivial. They follows immediately
from ones in theorems 2.1 and 3.1.\par

\bigskip

\section{Grand Lebesgue Spaces (GLS). }

\vspace{3mm}

 \ We recall here first of all for reader conventions some definitions and facts from
the theory of GLS spaces.\par

\vspace{3mm}

 \ Recently, see
\cite{Fiorenza1}, \cite{Fiorenza2},\cite{Ivaniec1}, \cite{Ivaniec2}, \cite{Jawerth1},
\cite{Karadzov1}, \cite{Kozatchenko1}, \cite{Liflyand1}, \cite{Ostrovsky1}, \cite{Ostrovsky2} etc.
 appear the so-called Grand Lebesgue Spaces (GLS)

  $$
G(\psi) = G(\psi,X) = G(\psi,X; A;B);  \ A;B = \const; \ A \ge 1, \ B \le \infty
 $$
spaces consisting on all the measurable functions $ f : X \to R  $ with finite norms

$$
||f||G(\psi) \stackrel{def}{=} \sup_{p \in (A;B)} \left[\frac{|f|_p}{\psi(p)} \right], \eqno(5.1)
$$
where as above

$$
|f|_{p,X} = |f|_p = \left[\int_X |f(x)|^p \ \mu(dx) \right]^{1/p}.
$$

 Here $ \psi = \psi(p), \ p \in (A,B) $ is some continuous positive on the {\it open} interval $ (A;B) $ function such
that

$$
\inf_{p \in(A;B)} \psi(p) > 0. \eqno(5.2)
$$

 \ We define formally  $ \psi(p) = +\infty $  for the values  $ \ p \notin [A,B]. $ \par

 \ We will denote also

$$
\supp(\psi) \stackrel{def}{=} (A;B).
$$

 \ The set of all such a functions with support $ \supp(\psi) = (A,B) $ will be denoted by  $ \ \Psi = \Psi(A,B). $  \par

 \ The  GLS space $ G(\zeta,Y) = G(\zeta,Y; A_1,B_1) $ based on the measurable space $ (Y,N,\nu) $ may be introduced quite
analogously. \par

 \ These spaces are complete and rearrangement invariant; and are used, for example, in
the theory of Probability, theory of Partial Differential Equations,
 Functional Analysis, theory of Fourier series, theory of
Martingales, Mathematical Statistics, theory of Approximation  etc. \par

 \ Notice that the classical Lebesgue-Riesz spaces $ L_r $  are extremal cases of the Grand Lebesgue Spaces. Indeed,
define the {\it degenerate} $ \Psi \ - $ function of the form $ \psi_{(r)}(p), \ r = \const \ge 1, \ p \ge 1 $ as follows:

$$
\psi_{(r)}(p):= 1, \ p = r; \hspace{4mm} \psi_{(r)}(p) = + \infty, \ p \ne r;
$$
and define formally $  C/\infty = 0, \ C = \const. $  Then

$$
||f||G(\psi_{(r)},X) =  | \ f  \ |_{r,X}.\eqno(5.3)
$$

\ The so - called {\it exponential} Orlicz spaces are also the particular cases of the GLS, see for instance
 \cite{Ostrovsky2},  \cite{Ostrovsky1}, chapter 1.  \par

 \ More detail, let for simplicity  $ \ \mu(X) = 1, \ $ i.e. let the measure $  \ \mu \ $ be probabilistic,
 and let the measurable function (random variable) $ \psi(\cdot) \in G\Psi = G\Psi_{\infty} $  be such that the new
 generated by $ \psi $ function

 $$
 \nu(p) = \nu_{\psi}(p) := p \ln \psi(p), \ p \in [1,\infty) \eqno(5.4)
 $$
is convex. The measurable function $ \eta = \eta(x) $ belongs to the space $  G\psi $ if and only if it  belongs to the Orlicz's
space $ L(N_{\psi}) $ with the correspondent {\it exponential} continuous Young - Orlicz function

$$
N_{\psi}(u) := \exp \left( - \nu^*_{\psi}(\ln |u|)  \right), \ |u| \ge e,  \eqno(5.5a)
$$

$$
N_{\psi}(u) := C \ u^2, \ |u| < e, \hspace{4mm} C \ e^2 = \exp \left( - \nu^*_{\psi}(1)  \right), \eqno(5.5b)
$$
and herewith of course both the Banach spaces norms: $  ||\cdot|| L(N_{\psi})  $  and $ || \cdot||G\psi  $ are equivalent. \par

 \ Here the $  \nu^*(\cdot) $ denotes as usually the Young-Fenchel, or Legendre transform of the function $ \nu(\cdot): $

$$
\nu^*(v) \stackrel{def}{=} \sup_{u > 0} (u \ v - \nu(u)), \ v > 0. \eqno(5.6)
$$

\vspace{3mm}

 \ One can also complete characterize (under formulated here conditions)  the belonging of the non - zero function $ f: X \to R $ to the
space $ G\psi $  by means of its tail behavior:

$$
f \in G\psi \Leftrightarrow \exists K = \const \in (0, \infty), \ \max( \mu \{x: f(x) > u \}, \mu \{x: f(x) < - u \}  ) \le
$$

$$
\exp \left( - \nu^*_{\psi}(\ln |u|/K)  \right), \ u \ge K e,
$$
see \cite{Kozatchenko1}, \cite{Ostrovsky1}, p. 33 - 35.\par

 \ For instance, the {\it random variable} $ \ \eta $ defined on some probability space $ \ (\Omega, B,{\bf P}), \ $
has a finite GLS norm of the form

$$
|| \ \eta \ ||_{m} \stackrel{def}{=} \sup_{p \ge 1} \left[ \frac{|\eta|_p}{p^{1/m}} \right] < \infty,
$$
where $ \ m \ $ is positive constant not necessary to be integer, if and only if

$$
\max( {\bf P}(\eta > u),  {\bf P}(\eta < - u)) \le  \exp\left(- C(m) \ u^m \right), \ u \ge 1.
$$

 \vspace{4mm}

 \  The case when  the supremum in (5.1) is calculated over {\it finite } interval is investigated in
 \cite{Liflyand1}, \cite{Ostrovsky8}:

$$
G_b \psi = \{\xi, \ ||\xi||G_b \psi < \infty \}, \  ||\xi||G_b \psi
\stackrel{def}{=} \sup_{ 1 \le p < b } \left[ \frac{|\xi|_p}{\psi(p)} \right], \  b = \const > 1,
$$
but here $ \psi = \psi(p) $ is  continuous  function in the semi-open interval $ 1 \le p < b $
such that $ \lim_{p \uparrow b} \psi(p) = \infty; $  the case when $ \psi(b - 0) < \infty $ is trivial. Indeed, if
$ \psi(b - 0) < \infty, $ then the space $ G_b \psi $ coincides up to norm equivalence with Lebesgue - Riesz  space $  \ L_b. $ \par

 \ We define formally in the case when $ b < \infty \ \psi(p) := + \infty \ $  for all the values $ p > b.$ \par

 \ Let a function $  f:  X \to R  $ be such that

 $$
 \exists (A,B): \ 1 \le A < B \le \infty, \ \forall p \in (A,B) \Rightarrow |f|_p < \infty.
 $$
 \ Then the function $  \psi = \psi(p) = \psi_f(p) $ may be {\it naturally} defined by the following way:

$$
\psi_f(p) := |f|_p, \ p \in (A,B). \eqno(5.7)
$$

 \ Evidently, $ \ ||\ f \ ||G\psi_f = 1. $ \par

\hfill $\Box$

\bigskip

\section{Acting of the composition operator on GLS. }

\vspace{3mm}

 \ {\it  Statement of problem. } Assume the function $  f(\cdot) $ belongs to some Grand Lebesgue pace $ G(\psi,Y) = G(\psi,Y; A_1,B_1), $ where
$  1 \le A_1 < B_1 \le \infty. $  Let also $  G(\zeta,X) = G(\zeta,X;A_2,B_2) $ be another GLS builded on the measurable space
$  (X,M,\mu). $ We set ourselves the problem of the norm estimate of product operator $ W_{g,\xi}[f] $ acting between two Grand
Lebesgue spaces

$$
D(\zeta, \psi) = D_{g, \xi}(\zeta, \psi)  \stackrel{def}{=} || \ W_{g, \xi} \ ||( \ G(\psi,Y)  \to G(\zeta,X) \ ). \eqno(6.1)
$$

 \ So, let $ \ f \in G(\psi,Y; A_1,B_1); \  $ we can and will suppose without loss of generality
 $ \ || \ f \ || G(\psi,Y; A_1,B_1) = 1.  $ This imply in particular

$$
\forall q \in (A_1, B_1) \ \Rightarrow |f|_{q,Y} \le \psi(q).
$$

 \ For instance, the function $ \ \psi = \psi(q) \ $  may be selected as a natural function for the function $  f(\cdot). $ \par

 \ We apply theorem 4.1:

$$
| \ W_{g, \xi} [f] \ |_{p,X} \le T_{g,\xi}(p,q) \ |f|_{q,Y} \le T_{g,\xi}(p,q) \ \psi(q). \eqno(6.2)
$$

 \ Introduce the following $ \ \Psi \- \  $  function

$$
\Theta(p) := \inf_{ q > p; q \in (A_1,B_1) } \left[   T_{g,\xi}(p,q) \ \psi(q) \right] \eqno(6.3)
$$
with correspondent support

$$
(A_2, B_2) := \supp \Theta(p), \eqno(6.4)
$$
then the   inequality (6.2) may be rewritten as follows

$$
| \ W_{g, \xi} [f] \ |_p \le \Theta(p), \ p \in (A_2, B_2). \eqno(6.5)
$$
 \ To summarize:

\vspace{3mm}

{\bf Theorem 6.1.}

$$
|| \ W_{g, \xi} \ ||(G(\psi,Y) \to G(\Theta,X) ) \le 1, \eqno(6.6)
$$
where the constant "1" in the right-hand side (6.6) is the best possible.\par

\vspace{3mm}

 \ The last assertion follows immediately from the main result of the article
 \cite{Ostrovsky100}, see also  \cite{Ostrovsky102}.\par

\vspace{4mm}

 \ {\bf Remark 6.1.} The multiplicative operators between two Orlicz's spaces are investigated in many works, see e.g.
\cite{Djakov1}, \cite{Lesnik1}, \cite{Maligranda1}, \cite{Shragin1}.\par

\bigskip

\section{Examples. }

\vspace{3mm}

 \ {\bf Example 7.1}. Consider the following product operator

$$
W^{(r,t)}[f](x) = x^{-1/t} \ f(x^r).  \eqno(7.1)
$$

 \ Here $  X = Y = (0,1), \ t = \const \in (0,1), r = \const > 0, r \ne 1. $\par

 \ We find using the proposition of theorem 4.1: $  \ || \ W^{(r,t)} ||(L_q \to L_p) = \phi(p,q;r,t),  $
where

$$
\phi(p,q;r,t) \stackrel{def}{=}
$$

$$
\inf_{ l \in ( \max(p,pr),q/(tq + 1)) } \left[ r^{-1/p} \times
\left\{ \ \left( \ \frac{l - p}{l- pr} \ \right)^{1/p - 1/l } \cdot \left( \ 1 - \frac{t l q}{q - l} \ \right)^{1/q - 1/l } \ \right\} \ \right]. \eqno(7.2)
$$

 \ Recall that here $  \ q > \max(p,pr). \ $ It will be presumed of course that $ \max(p,pr) < q/(tq + 1). $ \par

\vspace{3mm}

 \hspace{3mm} Let at first $  r < 1; $ we find then by some calculations

$$
\phi(p,q;r,t) \le r^{-1/p} \ e^{-1/q}.
$$

\vspace{3mm}

 \ The opposite case $ \ r > 1 \ $ is more complicated. Denote

$$
l_0 = l_0(p,q;r,t) := \frac{p^2 r(1 + qt) + q(r-1)}{(1 + q t)(r - 1 + tpr)}. \eqno(7.3)
$$

  \ It is easily to  verify that $ pr < l_0 < q/(qt + 1). $ Our statement:

$$
\phi(p,q;r,t) \le  r^{-1/p} \times
\left\{ \ \left( \ \frac{l_0 - p}{l_0 - pr} \ \right)^{1/p - 1/l_0} \cdot
\left( \ 1 - \frac{t l_0 q}{q - l_0} \ \right)^{1/q - 1/l_0 } \ \right\}.
$$

 \ Thus, the value $ l_0 = l_0(p,q;r,t) $ from (7.3) is asymptotically optimal in both the cases $ l \to pr + 0 $ and
 $ l \to q/(qt + 1) - 0. $  \par

\vspace{4mm}

 \ The expression for the function  $ \ \phi(p,q;r,t) \ $  allows a simplification. One can use the following identity

$$
\min_{x \in (a,b)} \left[ (x - a)^{-\gamma} \ (b - x)^{-\beta} \right] =
\frac{(\beta + \gamma)^{\beta + \gamma}}{\beta^{\beta} \ \gamma^{\gamma}} \cdot (b - a)^{ - (\beta + \gamma) }.
$$

 \ Here $ 0 < a  < b < \infty, \ \beta, \gamma = \const > 0.  $ \par

\vspace{3mm}

 \ We propose after some calculations in the case when $  q/(1 - qt) - pr \to 0, $ say
 $ \ 0 <  q/(1 - qt) - pr  < 1, \ $ and when $ 0 < t < 1, \ r > 1, $

$$
r_0 \le \min( p, q) \le \max(p,q) \le R_0, \ 0 < r_0 = \const < R_0 = \const < \infty:
$$

\vspace{3mm}

$$
\phi(p,q;r,t) \asymp \left\{ \frac{q}{1 - qt} - pr  \right\}^{-(1/q +t)}.
$$

\vspace{4mm}

 \ {\bf Example 7.2.A.} Multiplicative operator between Grand Lebesgue Spaces. \\

\vspace{3mm}

 \ Let now both the spaces $ \ (X,M,\mu) = (Y,N,\nu)  = ( [0,1], B, dx) \ $ be probability spaces. Assume for simplicity
  $ r = 1  $  and  $  t \in (0,1) $ in the examples (3.1), (and, after, (2.1)). Let also $ \ \psi = \psi(q), \ 1 \le q < \infty \ $ be certain
 $  \Psi \ - $  function with unbounded support. Introduce one still the following (linear) multiplicative operator acting on arbitrary
function $ \ f(\cdot) \ $ from the GLS space $ G\psi: \ f \in G\psi $

$$
V^{(t)}[f](y) \stackrel{def}{=} y^{-t} \cdot f(y).
$$

 \ Introduce for the values $ p \in [1, 1/t) $ the following $  \ \Psi \ - $ function

$$
\tau(p) := \inf_{ q > p/(1 - pt)} \left[ \ \left( 1- \frac{tpq}{q-p} \ \right)^{ (p-q)/pq  } \ \psi(q) \right]. \eqno(7.4)
$$

 \ We derive by virtue of the example 3.1 the following non - improvable in general case estimate of the form

$$
|| \ V^{(t)}[f] \ ||G\tau \le || \ f \ ||G\psi. \eqno(7.5)
$$

 \ Let us choose in (7.4)

 $$
  q = \frac{2 p}{1 - pt},
 $$
we find then

$$
\tau(p) \le \left[ \ \frac{1 + pt}{1 - pt} \ \right]^{ 2p/( 1 + pt ) }  \cdot \psi \left( \frac{2p}{1 - pt} \right). \eqno(7.6)
$$

 \ In particular, if $ \psi(q)  \asymp C_1 \ q^{1/m}, \  m = \const > 0, $ then

$$
\tau(p) \asymp C_2 (1/t - p)^{ - (1/t + 1/m) }, \ 1 \le p < 1/t. \eqno(7.7)
$$

 \ It is interest by our opinion to note that the function $ \tau(p) $ has a bounded support, despite the source $ \Psi \ -  $
function $ \psi(\cdot) $ has unbounded one.\par

\vspace{4mm}

 \ {\bf Example 7.2.B.} Composite operator between Grand Lebesgue Spaces. \\

\vspace{3mm}

\ Let again both the spaces $ \ (X,M,\mu) = (Y,N,\nu)  = ( [0,1], B, dx) \ $ be probability spaces. Assume for simplicity
  $ r > 1,  $ in the example (2.1), the so - called "hard case". \par

 \  Let also $ \ \psi = \psi(q), \ 1 \le q < \infty \ $ be certain
 $  \Psi \ - $  function. Consider as before  the following (linear) composite operator acting on arbitrary
function $ \ f(\cdot) \ $ from the GLS space $ G\psi: \ f \in G\psi $

$$
U_{(r)}[f](x):=  f(x^r).
$$

 \ Introduce for all the values $ p \in [1, \infty) $ the following $  \ \Psi \ - $ function

$$
\sigma(p)= \sigma_r[\psi](p):= r^{-1/p}  \inf_{ q > pr} \left\{ \left( \ \frac{q - p}{q - pr} \ \right)^{1/p \ - \ 1/q} \ \psi(q) \ \right\}.
 \eqno(7.8)
$$

 \ We derive by virtue of theorem 2.1 the following non-improvable in general case estimate of the form

$$
|| \ U_{(r)}[f] \ ||G\sigma \le || \ f \ ||G\psi. \eqno(7.9)
$$

 \ In particular, let us choose in (7.8) $  q:=q_0 = \lambda \ p \ r, \ \lambda = \const > r;  $ then

$$
\sigma_r[\psi](p) \le C_r(\lambda) \ \psi(\lambda \ p)
$$
or more precisely

$$
\sigma_r[\psi](p) \le r^{-1/p} \cdot \inf_{\lambda > r}
\left[ \ \left\{ \frac{\lambda - 1}{\lambda - r} \right\}^{(1/p)(1 - 1/r \lambda)} \cdot
 \ \psi(\lambda p) \ \right]. \eqno(7.10)
$$

 \ If for instance $  \psi(p) = \psi^{(m)}(p) = p^{1/m} $ or  more generally if the function $ \psi = \psi(p) $ satisfies the so - called
weak $ \ \Delta_2 \ $ condition at the infinity:

$$
\exists \ \lambda > r \ \Rightarrow  \sigma_r[\psi](p) \le C_{2,r}(\lambda) \ \psi(p),
$$
then

$$
|| \ U_{(r)}[\cdot] \ ||(G\psi \to G\psi)=  C_{3,r} < \infty.  \eqno(7.11)
$$

 \ On the other words, the function $  \ U_{(r)}[f](x)  $ may belong at the same space as the source function $ \ f(\cdot), \ $
in contradiction to the foregoing example. \par

 \ Notice that despite the function $ \ \psi(p) = \psi^{(m)}(p) = p^{1/m} \ $ satisfies the weak $ \ \Delta_2 \ $  condition
at the infinity, the correspondent Young - Orlicz function

$$
 N_{\psi^{(m)}}(u) \asymp \exp\left( C(m,r) \ |u|^m    \right), \ |u| \ge 1
$$
 does not. Since the Orlicz's norm $ \ L(  N_{\psi^{(m)}}), \ $ as we know, is equivalent to the $  G \psi^{(m)} $ norm, we conclude that

$$
|| \ U_{(r)}[\cdot] \ || \left\{L \left(N_{\psi^{(m)}} \right) \to L \left(N_{\psi^{(m)}} \right) \right\}  =  C\{4,r,m\} < \infty.
$$

 \ Note that this estimate does not follows from the main result of the article \cite{Cui1}. \par

\vspace{4mm}

 \ {\bf Example 7.3.} (Counterexample). \par

\vspace{3mm}

  Let again $ X = Y = (0,1)  $ and  define

 $$
 f(x) = x^{-1/2}, \ \hspace{6mm} \xi(x) = x^3. \eqno(7.12)
 $$
 Then

 $$
 |f|_p  = \left[ \frac{2}{2-p}   \right]^{1/p} =: \psi(p), \ 1 \le p < 2;
 $$

$$
z(x) = 3^{-1} x^{-2/3},  \ 0 < x \le 1; \hspace{5mm} |z|_q = 3^{1/p - 1} \cdot (3 - 2p)^{-1/p} =: \theta(p), \ 1 \le q < 3/2;
$$
or equally $ f(\cdot) \in G\psi, \ z(\cdot) \in G\theta,  $ but the superposition function $ g(x) = f(\xi(x)) = x^{-3/2}  $ does not
belongs to any $ L_p(X) $ space with $ p \ge 1.  $ \par
 \ The cause of seeming contradiction with theorem 2.1 is following: the function $ \ f(\cdot) \ $ does not belongs to arbitrary
Lebesgue - Riesz space $ \ L_q  \ $ with $ \ q > pr,\ p \ge 1, \ $  as long as here $ \ r = 3. $

\vspace{4mm}

{\bf Example 7. 4.} (Linear substituting). \par

\vspace{3mm}

 \ Let here $  X = R^d  $ with ordinary Lebesgue measure and $ f: R^d \to R  $ be some function  belonging to the space $  G \psi. $
 Let also $  A $ be non degenerate linear map (matrix) from $  R^d $ to itself. \par

 \ Define an operator of a view

 $$
 V_A[f] = f(Ax). \eqno(7.13)
 $$
 Obviously,

$$
| V_A[f]|_p^p = \int_{R^d} |f(A x)|^p \ dx = \int_{R^d}  |\det(A)|^{-1} \  |f(y)|^p \   dy =
$$

$$
|\det(A)|^{-1} \ |f|_p^p,
$$
or equally

$$
|V_A[f]|_p = |\det(A)|^{-1/p}  \ |f|_p,
$$
and following

$$
|V_A[f]|_p \le |\det(A)|^{-1/p} \cdot ||f||G\psi   \cdot \psi(p). \eqno(7.14)
$$

 \ Let the function $  \psi(\cdot) $ be factorable:

$$
\psi(p) = \frac{\zeta(p)}{ \tau(p)}, \ p \in (A,B),
$$
where both the functions  $ \zeta(\cdot), \ \tau(\cdot) $ are from the set $  G\Psi, $ i.e. satisfy all the conditions imposed
on the  function $  \psi(\cdot). $ We deduce after dividing the inequality (7.14) on the function $ \zeta(p): $

$$
\frac{|V_A[f]|_p}{\zeta(p)} \le ||f||G\psi \cdot \frac{[\det(A)]^{-1/p}}{\tau(p)}. \eqno(7.15)
$$

  \ Recall now that the fundamental  function  $ \phi(G\tau, \ \delta), \ 0 \le \delta \le \mu(X)  $
 for the Grand Lebesgue Space   $  G \tau  $ may be calculated by the formula

$$
\phi(G\tau, \ \delta) = \sup_{p \in (A,B)} \left[ \frac{\delta^{1/p}}{\tau(p)}  \right].
$$

 \ This notion play a very important role in the theory of operators, Fourier analysis etc., see \cite{Bennett1}.
The detail investigation of the fundamental  function  for GLS is done in \cite{Liflyand1},  \cite{Ostrovsky2}. \par

\vspace{3mm}

 Taking the maximum over $  p; \ p \in (A,B) $ from both the sides of inequality (7.15), we  get to the purpose of this
subsection:  under our condition

 $$
 ||V_A[f]|| G \zeta  \le ||f||G\psi \cdot \phi(G \tau, |\det(A)|^{-1}). \eqno(7.16)
 $$

\vspace{3mm}

 \ We note in conclusion that the multivariate case, for instance, the operator of the form

$$
W^{(r_1,r_2; \ t_1, t_2)}[f](x_1, \ x_2):= x_1^{-t_1} \ x_2^{-t_2} \ f(x_1^{r_1}, \ x_2^{r_2}),
$$
where $ \ r_1, r_2 = \const > 0, \  t_1, t_2 = \const \in (0,1), $ may be investigated analogously. \par


\bigskip

\begin{thebibliography}{99}

\vspace{4mm}

\bibitem{Abrahamse1}
{\sc M.B. Abrahamse.} {\it Mutiplication Operators.} Lecture Notes in Mathematics, 693, Springer
Verlag (1978), 17-36.

\bibitem{Anosov1}
{\sc D.V. Anosov.}  (2001) {\it Ergodic theory.} In Hazewinkel, Michiel, Encyclopedia of Mathematics, Springer, ISBN 978-1-55608-010-4.

\bibitem{Arora1}
{\sc S. C. Arora, Gopal Datt and Satish Verma.} {\it Composition operators on Lorentz spaces.} Bull.
Austral. Math. Soc. 76 (2007), 205-214.

\bibitem{Bennett1}
{\sc C. Bennett and R. Sharpley.} {\it Interpolation of operators.}
Orlando, Academic Press Inc., 1988.

\bibitem{Bedford1}
{\sc Tim Bedford, Michael Keane and Caroline Series, eds.} (1991). {\it Ergodic theory, symbolic dynamics and hyperbolic spaces.}
Oxford University Press, ISBN 0-19-853390-X,  (A survey of topics in ergodic theory; with exercises.)

\bibitem{Cui1}
{\sc Cui, Yunan; Hudzik, Henryk; Kumar, Romesh; Maligranda, Lech. } {\it Composition operators in
Orlicz spaces.} J. Aust. Math. Soc. 76 (2004), no. 2,  189-206.

\bibitem{Djakov1}
{\sc P. B. Djakov and M. S. Ramanujan.} {\it Multipliers between Orlicz sequence spaces.} Turk. J., Math., {\bf 24,}
(2000), 313-319.

\bibitem{Douglas1}
{\sc R.G. Douglas. } {\it Banach Algebra Techniques in Operator Theory. } Academic Press, New York,
1972.

\bibitem{Dunford1}
N.Dunford, B.Schwartz. {\it Linear Operators.} V.1, General Theory. Academic
Press, (1958), New York, London.

\bibitem{Estaremi1}
{\sc Y. Estaremi.  } {\it Multiplication and composition operators between two different Orlicz spaces.}\\
arXiv:1301.4830v1 [math.FA] 21 Jan 2013

\bibitem{Estaremi2}
{\sc Y. Estaremi and M.R. Jabbarzadeh. } {\it Weighted Lambert type operators on $ L_p \ - $ spaces. } Operators
and Matrices., {\bf 7,} (2013), 101-116.

\bibitem{Estaremi101}
{\sc Y.Estaremi, S.Maghsody and I.Rahmani.} {\it On properties of multiplication and composition operators between Orlicz spaces.}\\
arXiv:1506.00369v1 [math.FA] 1 Jun 2015

\bibitem{Fiorenza1}
{\sc A. Fiorenza.} {\it Duality and reflexivity in grand Lebesgue
spaces.} Collect. Math. {\bf 51}, (2000), 131-148.

\bibitem{Fiorenza2}
{\sc A. Fiorenza and G.E. Karadzhov. } {\it Grand and small Lebesgue
spaces and their analogs.} Consiglio Nationale Delle Ricerche,
Instituto per le Applicazioni del Calcoto Mauro Picone, Sezione
di Napoli, Rapporto tecnico,  272/03, (2005).

\bibitem{Giri1}
{\sc Ratan Kumar Giri, Shesadev Pradhan.}{\it Multiplication operators in Orlicz spaces. }
arXiv:1606.03363v1 [math.FA] 10 Jun 2016

\bibitem{Gupta1}
{\sc S. Gupta, B. S Komal and N. Suri.} {\it Weighted composition operators on Orlicz spaces.} Int. J.
Contemp. Math. Sciences, {\bf 1,} 11-20, \ (2010).

\bibitem{Ivaniec1}
{\sc T. Iwaniec and C. Sbordone.}  {\it On the integrability of the
 Jacobian under minimal hypotheses.} Arch. Rat. Mech. Anal., {\bf 119}, (1992), 129-143.

\bibitem{Ivaniec2}
{\sc T. Iwaniec, P. Koskela and J. Onninen. } {\it Mapping of Finite
Distortion: Monotonicity and Continuity.} Invent. Math. {\bf 144}, \ (2001), 507-531.

\bibitem{Jawerth1}
{\sc B. Jawerth and M. Milman.} {\it Extrapolation theory with
applications.} Mem. Amer. Math. Soc., {\bf 440}, \ (1991).

\bibitem{Karadzov1}
{\sc G.E. Karadzhov and M. Milman.} {\it Extrapolation theory: new
results and applications.} J. Approx. Theory, {\bf 113}, \ (2005),
38-99.

\bibitem{Kozatchenko1}
{\sc Yu.V. Kozatchenko and E.I. Ostrovsky.} {\it Banach spaces of random
variables of subgaussian type.} Theory Probab. Math. Stat., Kiev,
1985,  42-56, (in Russian).

\bibitem{Kumar1}
{\sc R. Kumar. } {\it Composition operators on Orlicz spaces.} Integral Equations Operator Theory, {\bf 29},
(1997), 17-22.

\bibitem{Lesnik1}
{\sc Karol Le\'snik and Jakub Tomaszewski.}  {\it  Pointwise multipliers of Orlicz function spaces and factorization. }\\
arXiv:1605.08581v1 [math.FA] 27 May 2016

\bibitem{Liflyand1}
{\sc Liflyand E., Ostrovsky E., Sirota L.} {\it Structural Properties of Bilateral Grand Lebesgue Spaces.}
Turk. J. Math.; {\bf 34}, (2010), 207-219.

\bibitem{Maligranda1}
{\sc L. Maligranda and E. Nakai.} {\it Pointwise multipliers of Orlicz spaces.} Arch. Math., {\bf 95,} (2010), no. 3,
251-256.

\bibitem{Montes-Rodrguez1}
{\sc Montes-Rodrguez, A.}  {\it The essential norm of a composition operator on Bloch spaces. } Pacific
J. Math., {\bf 188,} 339-351 (1999).

\bibitem{Musielak1}
{\sc J. Musielak.} {\it Orlicz spaces and modular spaces.} Lecture Notes in Math. 1034, Springer,
Berlin, (1983).

\bibitem{Okikiolu1}
{\sc G.O.Okikiolu.} {\it Aspects of the theory of bounded Integral Operators in the $ L^p $  Spaces. }
Academic Press; London,   New Yotk; (1971).

\bibitem{Ostrovsky1}
{\sc E.I. Ostrovsky.} {\it Exponential Estimations for Random Fields.}
Moscow-Obninsk, OINPE, 1999, (in Russian).

\bibitem{Ostrovsky2}
{\sc E. Ostrovsky and L.Sirota.} {\it Moment Banach spaces: theory and applications.}
HAIT Journal of Science and Engeneering, {\bf C}, Volume 4, Issues 1-2,
pp. 233-262, (2007).

\bibitem{Ostrovsky100}
{\sc E. Ostrovsky and L.Sirota.} {\it  Boundedness of operators in bilateral Grand Lebesgue Spaces,
with exact and weakly exact constant calculation. } \\
arXiv:1104.2963v1 [math.FA] 15 Apr 2011

\bibitem{Ostrovsky101}
{\sc E. Ostrovsky and L.Sirota.}
{\it  Multiple weight Riesz and Fourier transforms  in bilateral anosotropic Grand
Lebesgue Spaces.  } \\
arXiv:1208.2392v1 [math.FA] 12 Aug 2012

\bibitem{Ostrovsky102}
{\sc E.Ostrovsky, L.Sirota, E.Rogover.} {\it Integral Operators in Bilateral Grand Lebesgue Spaces.} \\
arXiv:0912.2538v1 [math.FA]

\bibitem{Ostrovsky3}
{\sc E.Ostrovsky,  L.Sirota, E.Rogover.} {\it  Riesz's and Bessel's operators in Bilateral Grand Lebesgue Spaces. } \\
arXiv:0907.3321v1 [math.FA] 19 Jul 2009

\bibitem{Ostrovsky4}
{\sc E.Ostrovsky,  L.Sirota.} {\it Weight Hardy - Littlewood inequalities for different powers. } \\
arXiv:0910.5880v1 [math.FA] 30 Oct 2009

\bibitem{Ostrovsky6}
{\sc E.Ostrovsky,  L.Sirota.} {\it  Compact sets in bide-side Grand Lebesgue Spaces, with applications.}
arXiv:0902.2916v1 [math.FA] 17 Feb 2009

\bibitem{Ostrovsky7}
{\sc  Ostrovsky E. and Sirota L.}
{\it Moment Banach spaces: theory and applications.}
HIAT Journal of Science and Engineering, {\bf C}, Volume 4, Issues 1 - 2,
pp. 233-262, (2007).

\bibitem{Ostrovsky8}
{\sc E.Ostrovsky,  L.Sirota, E.Rogover.} {\it  Integral operators in bilateral Grand Lebesgue Spaces.} \\
arXiv:0912.2538v1 [math.FA] 13 Dec 2009

\bibitem{Petersen1}
{\sc Karl Petersen.} {\it Ergodic Theory}. \ Cambridge Studies in Advanced Mathematics. Cambridge: Cambridge University Press., 1990.

\bibitem{Rao1}
{\sc Rao M.M., Ren Z.D.} {\it Theory of Orlicz Spaces.} Marcel Dekker Inc., 1991. New
York, Basel, Hong Kong.

\bibitem{Shragin1}
{\sc I. V. Shragin. }  {\it On certain operators in generalized Orlicz spaces.} Dokl. Akad. Nauk SSSR \ (N.S:), {\bf 117,}
(1957), 40-43, (in Russian).

\bibitem{Shapiro1}
{\sc J. H. Shapiro.} {\it The essential norm of a composition operator.} Analls of Math., {\bf 125,} \ (1987),
375-404.

\bibitem{Singh1}
{\sc R. K. Singh and J. S. Manhas. } {\it Composition operators on function spaces.} North Holland
Math. Studies, {\bf 179,} \ Amsterdam, \ 1993.

\bibitem{Stein1}
{\sc E.M.Stein.} {\it Singular Integrals and Differentiability Properties  of Functions. }
Princeton University Press, Princeton, \ (1992).

\bibitem{Takagi1}
{\sc  H. Takagi. } {\it Compact weighted composition operators on } $ L_p.$  Proc. Amer. Math. Soc., {\bf 116,}
(1992), 505-511.

\bibitem{Takagi2}
{\sc H.Takagi.} {\it Fredholm Weighted Composition Operators.} Integral Equations and Operator
Theory, {\bf 16,} (1993), 267-276.

\bibitem{Takagi3}
{\sc H.Takagi, T., Miura, T. and Takahasi, Sin-Ei.} {\it Essential norm and stability constants of weighted
composition operators on C(X). } Bull. Korean Math. Soc., {\bf 40,} \  583-591, (2003).

\bibitem{Takagi4}
{\sc H.Takagi and K. Yokouchi.} {\it Multiplication and composition operators between two $ L_p \ - $ spaces. }
Contemporary Math. {\bf 232, }(1999), 321-338.

\bibitem{Walters1}
{\sc Walters, Peter.} (1982). {\it An introduction to ergodic theory. } Graduate Texts in Mathematics, {\bf 79, \ } Springer-Verlag,
ISBN 0-387-95152-0, Zbl. 0475.28009

\bibitem{Zheng1}
{\sc Zheng, L. } {\it  The essential norms and spectra of composition operators on}  $ H_1. $ Pacific J. Math., \
{\bf 203,} \  503-510, \ (2002).

\end{thebibliography}
\end{document}